\documentclass[11pt]{article}

\usepackage{amscd,amsmath, amssymb, fancyhdr, epsfig, color, tocloft, mathtools}
\usepackage{dutchcal}

\usepackage[backref=page]{hyperref}
\renewcommand*{\backrefalt}[4]{%
	\ifcase #1 (Not cited.)%
	\or        (Cited on page~#2.)%
	\else      (Cited on pages~#2.)%
	\fi}

\hypersetup{
	colorlinks   = true,
	citecolor    = magenta,
	linkcolor    = blue,
	urlcolor     = magenta	
}

\newcommand{\version}{version 2.0,\ \ Mar. 28, 2023}

\makeatletter
\def\x@arrow{\DOTSB\Relbar}
\def\xlongequalsignfill@{\arrowfill@\x@arrow\Relbar\x@arrow}
\providecommand{\xlongequal}[2][]{%
	\ext@arrow 0099\xlongequalsignfill@{#1}{#2}}
\def\xlongrightarrowfill@{\arrowfill@\relbar\relbar\longrightarrow}
\makeatother

\numberwithin{equation}{section}

\newcommand{\GL}{{\sf{GL}}}

\newcommand{\U}{{\sf{U}}}

\def\eqref#1{(\ref{#1})}

\newcommand{\goth}{\mathfrak}

\newcommand{\N}{{\mathbb N}}
\newcommand{\Z}{{\mathbb Z}}
\newcommand{\C}{{\mathbb C}}
\newcommand{\R}{{\mathbb R}}

\def\1{\sqrt{-1}\,}

\newcommand{\restrict}[1]{{\left|_{{\phantom{|}\!\!}_{#1}}\right.}}

\newcommand{\cntrct}                
{\hspace{2pt}\raisebox{1pt}{\text{$\lrcorner$}}\hspace{2pt}}

\newcommand{\arrow}{{\:\longrightarrow\:}}

\newcommand{\calo}{{\cal O}}

\renewcommand{\bar}{\overline}
\renewcommand{\phi}{\varphi}
\renewcommand{\epsilon}{\varepsilon}
\renewcommand{\geq}{\geqslant}
\renewcommand{\leq}{\leqslant}

\newcommand{\im}{\operatorname{im}}
\newcommand{\End}{\operatorname{End}}

\newcommand{\Tot}{\operatorname{Tot}}
\newcommand{\Id}{\operatorname{Id}}

\renewcommand{\sup}{{\operatorname{{\sf sup}}}}

\newcommand{\Aut}{\operatorname{Aut}}

\renewcommand{\dim}{\operatorname{\sf dim}}

\newcommand{\Spec}{\operatorname{{\sf Spec}}}

\newcounter{Mycounter}[section]
\newcounter{lemma}[section]
\renewcommand{\thelemma}{{Lemma \thesection.\arabic{lemma}}}
\newcommand{\lemma}{%
	\setcounter{lemma}{\value{Mycounter}}
	\refstepcounter{lemma}
	\stepcounter{Mycounter}
	{\noindent \bf \thelemma:\ }}

\newcounter{claim}[section]
\renewcommand{\theclaim}{{Claim \thesection.\arabic{claim}}}
\newcommand{\claim}{%
	\setcounter{claim}{\value{Mycounter}}
	\refstepcounter{claim}
	\stepcounter{Mycounter}
	{\noindent \bf \theclaim:\ }}

\newcounter{sublemma}[section]
\newcounter{corollary}[section]
\renewcommand{\thecorollary}{{Corollary \thesection.\arabic{corollary}}}
\newcommand{\corollary}{%
	\setcounter{corollary}{\value{Mycounter}}
	\refstepcounter{corollary}
	\stepcounter{Mycounter}
	{\noindent \bf \thecorollary:\ }}

\newcounter{theorem}[section]
\renewcommand{\thetheorem}{{Theorem \thesection.\arabic{theorem}}}
\newcommand{\theorem}{%
	\setcounter{theorem}{\value{Mycounter}}
	\refstepcounter{theorem}
	\stepcounter{Mycounter}
	{\noindent \bf \thetheorem:\ }}

\newcounter{conjecture}[section]
\newcounter{proposition}[section]
\renewcommand{\theproposition} {{Proposition \thesection.\arabic{proposition}}}
\newcommand{\proposition}{%
	\setcounter{proposition}{\value{Mycounter}}
	\refstepcounter{proposition}
	\stepcounter{Mycounter}
	{\noindent \bf \theproposition:\ }}

\newcounter{definition}[section]
\renewcommand{\thedefinition} {{Definition~\thesection.\arabic{definition}}}
\newcommand{\definition}{%
	\setcounter{definition}{\value{Mycounter}}
	\refstepcounter{definition}
	\stepcounter{Mycounter}
	{\noindent \bf \thedefinition:\ }}

\newcounter{example}[section]
\renewcommand{\theexample}{{Example \thesection.\arabic{example}}}
\newcommand{\example}{%
	\setcounter{example}{\value{Mycounter}}
	\refstepcounter{example}
	\stepcounter{Mycounter}
	{\noindent \bf \theexample:\ }}

\newcounter{remark}[section]
\renewcommand{\theremark}{{Remark \thesection.\arabic{remark}}}
\newcommand{\remark}{%
	\setcounter{remark}{\value{Mycounter}}
	\refstepcounter{remark}
	\stepcounter{Mycounter}
	{\noindent \bf \theremark:\ }}

\newcounter{problem}[section]
\newcounter{question}[section]
\makeatletter

\def\blacksquare{\hbox{\vrule width 5pt height 5pt depth 0pt}}
\def\endproof{\blacksquare}

\newcommand{\proof}{{\bf Proof: \ }}
\newcommand{\pstep}{{\bf Proof. Step 1: \ }}

\begin{document}

\begin{center}
{\Large\bf Algebraic cones of LCK manifolds with potential}\\[5mm]
{\large
Liviu Ornea\footnote{Liviu Ornea is  partially supported by Romanian Ministry of Education and Research, Program PN-III, Project number PN-III-P4-ID-PCE-2020-0025, Contract  30/04.02.2021},  
Misha Verbitsky\footnote{Misha Verbitsky is partially supported by
 the HSE University Basic Research Program, FAPERJ E-26/202.912/2018 
and CNPq - Process 313608/2017-2.\\[1mm]

}\\[4mm]

}

\end{center}

{\small
\hspace{0.15\linewidth}
\begin{minipage}[t]{0.8\linewidth}
{\bf Abstract} \\ 
A complex manifold $X$, $\dim X >2$, 
is called ``an LCK manifold with potential'',
if it can be realized as a complex submanifold of a Hopf manifold.
Let $\tilde X$ be its $\Z$-covering, considered
as a complex submanifold in $\C^n \backslash 0$.
We prove that $\tilde X$ is algebraic. 
We call the manifolds obtained this way {\bf the algebraic cones},
and show that the affine algebraic structure on $\tilde X$ is independent
from the choice of $X$. We give several intrinsic definitions
of an algebraic cone, and prove that these definitions are
equivalent.
\end{minipage}
}

\tableofcontents

\section{Introduction}

\subsection{LCK manifolds with potential}

An LCK manifold with potential (\cite{ov_lckpot}) can be defined in terms
of a certain Hermitian metric (\ref{_LCK_pote_def_}).
For our purpose, another definition is more convenient.
In complex dimension $>2$, an LCK manifold with potential 
can be defined as a compact complex manifold which admits a holomorphic
embedding to a linear Hopf manifold 
(\ref{_linear_Hopf_Definition_}).\footnote{The
equivalence of these two definitions 
is proven in dimension $>2$ (\cite{ov_indam}); for complex surfaces
it holds up to Kato conjecture, which is one of the greatest
unresolved conjectures in the geometry of complex surfaces.
Note that the Kato conjecture is widely believed to be true.}

Let $M$ be an LCK manifold with potential,
embedded to a Hopf manifold $\frac{\C^n \backslash 0}{\langle A \rangle}$.
Then the $\Z$-cover $\tilde M$ of $M$
is a complex submanifold in $\C^n \backslash 0$.

The closure $\tilde M_c$ of $\tilde M\subset \C^n \backslash 0$ in $\C^n$
is a Stein variety with an isolated singularity, called
{\bf the (weak) Stein completion of $\tilde M$} 
(\ref{_Stein_completion_Definition_}). The nature
of this singularity was not explored so far; in the
present paper, we relate its geometry to the classical
subject of projective normality, due to Zariski, Muhly
and Hodge-Pedoe (\cite{_Muhly_,_Zariski:complete_,_Hodge_Pedoe_}). 

The purpose of this paper is to study the algebraic geometry
of the variety $\tilde M_c$. We use the LCK geometry mostly as a motivation;
almost all the results are based on  arguments of complex-analytic
and algebraic geometric nature which do not invoke the 
differential geometry of an LCK manifold. 

\subsection{Algebraic cones}

Locally conformally K\"ahler manifolds
with potential are in a sense very ``algebraic''.
The difference between complex projective and K\"ahler geometry
is understood through the Kodaira embedding theorem:
a K\"ahler manifold with rational K\"ahler class
admits a complex embedding into a projective space.
In this setup, one can say that the geometry of
complex projective manifolds is part of the K\"ahler geometry,
and this part is certainly very algebraic.

The analogue of Kodaira embedding theorem is already known:
an LCK manifold with potential is precisely a complex manifold
which admits a holomorphic embedding to a Hopf manifold.
This notion is already very close to algebraic geometry.
For example, there exists a $p$-adic version of the theory
of Hopf manifolds and their complex subvarieties,
in the framework of the rigid analytic spaces
(\cite{_Scholze:congress_}). The $p$-adic Hopf manifolds are known 
and well-studied (\cite{_Mustafin_,_Voskuil_}).

The passage from LCK geometry to complex algebraic geometry
is based on the notion of the ``algebraic cone'', introduced
in \cite{ov_pams}. The {\bf closed algebraic cone} of an LCK manifold
with potential is its (weak) Stein completion $\tilde M_c$
(\ref{_Stein_completion_Definition_}),
equipped with a structure of an affine variety,
and {\bf the open algebraic cone} is $\tilde M$,
that is, $\tilde M_c$ without its origin (apex) point.
We give an intrinsic definition of the closed and
open algebraic cones in Subsection \ref{_cone_Subsection_}.

\subsection{Projective normality}

Let $X$ be a projective variety, and
$ R:= \bigoplus_{i\in \Z^{\geq 0}} H^0(X, \calo(i))$ 
the corresponding graded ring. Since $R$ is finitely
generated, its spectrum is an affine variety,
called {\bf the affine cone} of $X$.

Recall that a complex variety $X$ is called 
{\bf normal} if any bimeromorphic finite 
map $X' \arrow X$ is an isomorphism.\footnote{We 
actually use an equivalent definition (Subsection
\ref{_normality_of_cone_Subsection_}): a variety is normal if all
locally bounded meromorphic functions
are holomorphic.}
We discuss the normality in some detail in 
Section \ref{_normality_of_cone_Subsection_}.

This is related to the notion of ``projective normality''
much studied in classical algebraic geometry.
A projective variety $X \subset \C P^n$ is called 
{\bf projectively normal} if the ring
$R= \bigoplus_{i\in \Z^{\geq 0}} H^0(X, \calo(i))$ of the
regular functions in its affine cone is integrally closed;
this is equivalent to the affine cone being 
normal. The projective normality is a very subtle notion,
because it greatly depends on the choice of
projective embedding (\cite[ Exercise I.3.18]{_Hartshorne:AG_}).
Before the advent of the Hodge theory,
projective normality was used to prove things such as the Riemann-Roch formula.
O. Zariski and his student H. T. Muhly, who were the first to study
the projective normality in a systematic manner,
called it ``arithmetic normality'' 
(\cite{_Muhly_,_Zariski:complete_});
however, the ``projective normality'' stuck 
because it was used in the textbook by Hodge and
Pedoe, \cite{_Hodge_Pedoe_}. 

\subsection{New results}

In this subsection, we summarize the new results of this paper.

From the Remmert-Stein theorem (\cite[Chapter 2, \S 8.2]{demailly})
it follows that the closure $\tilde M_c$ of $\tilde M$ in $\C^n$ 
is a Stein variety with a singular point in the origin 
(\ref{_Stein_Remmert_completion_Remark_}).
We call this closure {\bf a weak Stein completion} of $\tilde M$.
Note that a weak Stein completion Ăâem a prioriÎ depends on the
embedding $M \hookrightarrow H$;
however, for an appropriate embedding, the completion
$\tilde M_c$ is normal, and this normal variety,
called {\bf the Stein completion},
is uniquely determined by the complex geometry of 
$\tilde M$ (\ref{_Forster_Stein_uniqueness_Remark_}). 

We prove that the Stein completion $\tilde M_c$
homeomorphically and holomorphically
projects to all weak Stein completions
(\ref{_Stein_completion_second_proof_Proposition_}).

We consider $\tilde M_c$ as a Stein variety with a
unique isolated singular point. 
In \cite{ov_lckpot,ov_pams} we proved that this Stein variety
is in fact an affine cone over a projective
variety (\ref{potcon}). In this paper we describe this cone
from the algebraic-geometric point of view.

In \ref{_algebra_structure_unique_Theorem_}, we prove that 
the (affine) algebraic structure of $\tilde M$
is uniquely determined by its complex analytic structure,
as long as the $\Z$-action is algebraic. 
{\it Ipso facto} this algebraic structure is independent
from the choice of the LCK manifold with potential
$M= \frac{\tilde M}{\Z}$ used to construct $\tilde M$.

We define {\bf a closed algebraic cone}
as an affine cone over a projective orbifold,
which is smooth outside of the origin. 
The variety $\tilde M_c$ is an example of a closed
algebraic cone.
In \ref{_cones_into_Hopf_equiv_Theorem_} we prove that 
this definition is equivalent to the definition of the
closed algebraic cone given in \cite{ov_pams} 
in terms of LCK manifolds with potential.

In this paper we study the set of closed
algebraic cones which correspond to the same
open algebraic cone.

We use the Stein completion to prove that
any open algebraic cone $\tilde M$ can be associated
with a unique closed algebraic cone $\tilde M_c$, which is
normal (\ref{_Stein_completion_second_proof_Proposition_}). 
We also prove that the other closed algebraic cones associated
to $\tilde M$ are homeomorphic to $\tilde M_c$,
because they are all obtained by adding a point
in the origin. We prove that the natural homeomorphism map 
$n:\; \tilde M_c\arrow \tilde M_c'$ to another closed 
algebraic cone associated to $\tilde M$
is holomorphic. Moreover, $n$ is the
normalization map.

\section{Closed and open algebraic cones}
\label{_cones_Section_}

In this section we give several definitions of
closed and open algebraic cones. 
In the main body of this paper we 
prove that all these definitions are equivalent.

\subsection{Stein completions}

Let $K\subset \C^n$ be a compact subset.
Hartogs theorem implies that any holomorphic function
on $\C^n \backslash K$ can be extended to $\C^n$, if $n
\geq 2$. Due to Rossi, the same
result is true for any Stein variety:

\hfill

\theorem\label{_Hartogs_Stein_Theorem_}
Let $X$ be a normal Stein variety, $\dim_\C X >1$, 
and $K \subset X$ a compact subset. Then every holomorphic
function on $X \backslash K$ can be extended to $X$.

\proof \cite[Theorem 6.6]{_Rossi:fields_}. \endproof

%

\hfill

\remark \label{_Forster_Stein_uniqueness_Remark_}
Let $A$ be a commutative Fr\'echet algebra 
over $\C$ and $\Spec(A)$
the {\bf continuous spectrum} of $A$, defined as the
set of all continuous $\C$-linear homomorphisms $A\arrow \C$.
By \cite{_Forster:1966_,_Forster:1967_},
$\Spec(H^0(\calo_X))=X$ for any Stein
variety $X$, where $H^0(\calo_X)$ is the algebra
of holomorphic functions equipped with the topology
of uniform convergence on compacts.
From \cite{_Forster:1967_} it
also follows that a Stein manifold is determined uniquely
from  $H^0(\calo_X)$ considered as a ring with $C^0$-topology.

\hfill

\definition\label{_Stein_completion_Definition_}
Let $X$ be a normal Stein variety, and $K\subset X$ a compact subset.
By \ref{_Hartogs_Stein_Theorem_}, the ring of functions on 
$X$ is identified with $H^0(\calo_{X \backslash K})$;
by \ref{_Forster_Stein_uniqueness_Remark_}, 
this ring with $C^0$ topology uniquely defines $X$.
Following \cite{andreotti_siu},
we call $X$ {\bf the Stein completion} of
$X\backslash K$. {\bf A weak Stein completion}
of $X \backslash K$ is any Stein variety
$X'$ containing a compact $K'$ such that 
$X'\backslash K'\cong X\backslash K$.

\hfill

\remark\label{_normalization_of_weak_Stein_comple_Remark_}
Clearly, if $X\backslash K$ is normal,
the normalization of a weak Stein completion
of $X \backslash K$ is the Stein completion
of  $X \backslash K$. 

\hfill

\example
Let $M \subset H$ be a submanifold in a Hopf manifold
$H= \frac{\C^N \backslash 0}{\Z}$, and
$\tilde M\subset \C^N \backslash 0$ its $\Z$-cover.
By the Remmert-Stein theorem 
(\ref{_Stein_Remmert_completion_Remark_}), 
the closure of $\tilde M$ is Stein;
however, it is not necessarily normal 
(\ref{_non_unique_non_normal_Example_}).
This is why we need the notion of a ``weak Stein 
completion''.

\hfill

\example\label{_Stein_comple_cone_Example_}
Let $P$ be a projective orbifold, and
$L$ an ample line bundle on $P$.
Consider the total space  $\Tot^\circ(L)$
of all non-zero vectors in $L$. 
By \ref{_Stein_completion_cone_Proposition_} and 
\ref{_Stein_completion_second_proof_Proposition_}
below, its (weak) Stein completion is obtained by 
adding a single point, called {\bf the apex},
or {\bf the origin}. The weak Stein completion
of $\Tot^\circ(L)$ is not unique, as 
\ref{_non_unique_non_normal_Example_} shows;
however, all weak Stein completions of $\Tot^\circ(L)$
are homeomorphic (\ref{_Stein_completion_second_proof_Proposition_}).

\subsection{Algebraic structures on algebraic cones}
\label{_cone_Subsection_}

The notion of an algebraic cone has two flavours:
there are ``closed algebraic cones'', which can be obtained
as weak Stein completions of ``open algebraic cones''.
The latter can be defined in algebraic terms as follows.

\hfill

\definition
Let $P$ be a projective orbifold, and
$L$ an ample line bundle on $P$.
Assume that the total space $\Tot^\circ(L)$
of all non-zero vectors in $L$ is smooth.
{\bf An open algebraic cone} is $\Tot^\circ(L)$.
The corresponding {\bf closed algebraic cone}
is its weak Stein completion $Z$, such that the
natural $\C^*$-action can be extended to $Z$
holomorphically. By \ref{_Stein_completion_second_proof_Proposition_}, the
closed algebraic cone is obtained by adding 
one point, called ``the apex'', or ``the origin''.

\hfill

%

Further on, we shall rely on the notion of ``holomorphic
contraction''. To simplify the terminology, 
we tacitly assume that the ``contractions'' are 
invertible; however, most of the arguments remain valid
without this assumption.

\hfill

\definition 
A {\bf contraction} of a topological space $X$ to a point
$x\in X$ is a homeomorphism $\phi:\; X \arrow X$ such that for any 
compact subset $K\subset X$ and any open set $U\ni x$, 
there exists $N>0$ such that for all
$n>N$, the map $\phi^n$ maps $K$ to $U$.

\hfill

{\bf An algebraic structure} on a 
complex analytic variety $Z$ is a subsheaf of the
sheaf of holomorphic functions which can be
realized as a sheaf of regular functions 
for some biholomorphism between $Z$ and
a  quasi-projective variety.

\hfill

One of the main results of this paper is the
following theorem, proven in Section \ref{_algebraic_comple_Section_}
below. Note that any closed algebraic cone $X$
admits an algebraic contraction 
$\gamma:\; X \arrow X$ taking $v$ to
$\lambda v$, for $\lambda$ an invertible
complex number, $|\lambda|<1$.

\hfill

\theorem\label{_algebra_structure_unique_Theorem_}
Let $X$ be a closed algebraic cone,
and ${\cal A}_1$, ${\cal A}_2$
two algebraic structures on $X$.
Assume that the algebraic varieties
$(X, {\cal A}_1)$ and $(X, {\cal A}_2)$
are affine and admit (invertible) algebraic contractions 
$\gamma_1$ and $\gamma_2$,
that is, holomorphic contractions 
compatible with the respective algebraic
structures. Then ${\cal A}_1={\cal A}_2$.

\hfill

{\bf Proof:} By \ref{_cones_into_Hopf_equiv_Theorem_}
below, $X$ can be
obtained in two different
ways as a weak Stein completion of $\Z$-coverings
of LCK manifolds with potential associated
with the holomorphic contractions $\gamma_1$ and $\gamma_2$.
Then  ${\cal A}_1={\cal A}_2$ by
\ref{_same_algebra_structure_on_cone_Theorem_}. \endproof

\subsection{Algebraic cones and subvarieties in Hopf manifolds}

The algebraic cones can be defined in terms
of projective orbifolds as above; there are two
more definitions equivalent to this one. We could
define open algebraic cones as $\Z$-coverings of submanifolds of a 
linear Hopf manifold (which is done in the present
section) or define the closed algebraic cones
as Stein varieties with an isolated
singularity admitting a holomorphic contraction, as in 
Subsection \ref{_cones_contraction_Subsection_} below.

\hfill

\definition\label{_linear_Hopf_Definition_}
A {\bf linear Hopf manifold} is a complex manifold
$H:= \frac{\C^n \backslash 0} {\langle A \rangle}$,
where $A\in \GL(n, \C)$ is an invertible linear contraction.
A {\bf diagonal Hopf manifold} is a linear
Hopf manifold such that $A$ is diagonalizable,
and {\bf a classical Hopf manifold} is 
 a linear Hopf manifold such that $A$ is
a scalar matrix.

\hfill

\theorem\label{_cones_into_Hopf_equiv_Theorem_}
Let $M \subset H$ be a submanifold in a Hopf manifold,
and $\tilde M\subset \C^n \backslash 0$ its $\Z$-covering.
Then $\tilde M$ is an open algebraic cone. Moreover,
any open algebraic cone can be obtained this way.

\hfill

\proof
Let $\tilde M$ be an open algebraic cone
equipped with the standard $\C^*$-action, and 
$\gamma\in \Aut(\tilde M)$ the automorphism
associated with $\lambda \in \C^*$, $|\lambda| <1$.
By \ref{_alge_cone_to_Hopf_Theorem_},
$\frac{\tilde M}{\langle \gamma\rangle}$
admits an embedding to a linear Hopf manifold.
To obtain the converse assertion, we 
use the locally conformally K\"ahler metrics
with potential, defined below (\ref{_LCK_pote_def_}).
As shown in \cite{ov_pams}, any submanifold
of a linear Hopf manifold is LCK with (proper) potential.
\ref{_cone_cover_for_LCK_pot_Theorem_}
implies that a K\"ahler $\Z$-cover of 
an LCK manifold with proper potential is
an open algebraic cone.
\endproof

\subsection{Algebraic cones and holomorphic contractions}
\label{_cones_contraction_Subsection_}

Another equivalent definition of algebraic cones
uses holomorphic contractions. Note that this
definition does not refer to algebraic structures,
however, an algebraic structure
compatible with a contraction 
is unique by \ref{_algebra_structure_unique_Theorem_}.

\hfill

\theorem
Let $X$ be a Stein variety with a unique isolated singularity,
admitting a holomorphic contraction. Then $X$ is
biholomorphic to a closed algebraic cone. Conversely,
any closed algebraic cone can be obtained this way.

\hfill

\proof
By the same argument as used in the proof of 
\ref{_Stein_completion_second_proof_Proposition_},
a closed algebraic cone $X$ can be obtained as the closure 
of the preimage $\pi^{-1}(P)$,
where $P \subset \C P^n$ is a projective orbifold
and $\pi:\; \C^{n+1}\backslash 0 \arrow \C P^n$
the standard projection. Then $X$ admits
a contraction $v \arrow \lambda v$, where $|\lambda|<1$.

Conversely, suppose that $X$ is 
a Stein variety  with a unique isolated singularity $x$ 
admitting a holomorphic contraction $\gamma$ to $x$.
Then 
$\frac{X \backslash \{x\}}{\langle \gamma\rangle}$
is a complex manifold admitting a holomorphic embedding 
to a Hopf manifold (\ref{_alge_cone_to_Hopf_Theorem_}), hence it is an 
LCK manifold with proper potential (\ref{potcon}).
By \ref{_cone_cover_for_LCK_pot_Theorem_}, 
$X\backslash \{x\}$
is an open algebraic cone.
\endproof

\section{Locally conformally K\"ahler manifolds}

The main motivation for this paper comes from the
theory of LCK (locally conformally K\"ahler) 
manifolds, though we could state and
prove everything in it without using this important
notion.

\hfill


\definition
A complex  manifold $(M, I)$ is called {\bf locally conformally
	K\"ahler} (LCK, for short) if it admits a
covering $(\tilde M, I)$ equipped with a K\"ahler metric
$\tilde\omega$ such that the deck group of the
cover acts on $(\tilde M, \tilde \omega)$ by holomorphic 
homotheties. An {\bf LCK metric} on an LCK manifold
is an Hermitian metric on $(M,I)$ such that its
pullback to $\tilde M$ is conformal with $\tilde \omega$.

\hfill

\remark Equivalently, a Hermitian manifold $(M,I,g)$ is
LCK if there exists a closed 1-form $\theta$ (called {\bf
  the Lee form}) such that the fundamental form
$\omega(x,y):=g(Ix,y)$ satisfies
$d\omega=\theta\wedge\omega$.

\hfill

The special class of LCK manifolds we are concerned with is the following:

\hfill

\definition\label{_LCK_pote_def_} 
An LCK manifold has
                 {\bf a proper LCK potential} if it
admits a K\"ahler $\Z$-covering on which the K\"ahler metric
has a global and positive  potential function $\psi$ such that
the deck group multiplies $\psi$ by a
constant.\footnote{Such a function is called {\bf automorphic}.}
In this case, $M$ is called {\bf an LCK manifold with
proper	potential}.

\hfill

In the sequel, we will often tacitly omit the word ``proper''.

\hfill

\example \label{_LCK_pot_Example_}
The linear Hopf manifolds are LCK with potential,
\cite{ov_non_linear}. For a classical Hopf manifold
$H:=(\C^n\setminus 0)/\langle A\rangle$, $A=\lambda\Id$, 
$|\lambda|> 1$,  the flat K\"ahler metric $\tilde g_0=\sum
dz_i\otimes d\bar z_i$ on $\C^n$ is multiplied by $\lambda^2$ by the
deck group $\Z$. Also, $\tilde g_0$ has the global
automorphic potential $\psi:=\sum |z_i|^2$.

%
%

%
%
%
%

\hfill

The main properties of LCK manifolds with 
(proper) potential are the following.

\hfill

\theorem \label{potcon}
{(\cite{ov_lckpot,ov_pams})}
Let $M$ be a compact LCK manifold with proper potential,
and $\tilde M$ its K\"ahler $\Z$-cover.
If  $\dim_\C M\geq 3$, then the metric 
completion $\tilde M_c$
is identified with the Stein completion of $\tilde M$,
the complement $\tilde M_c\setminus  \tilde M$
is just one point. Moreover, $\tilde M_c$ is an affine
algebraic variety obtained as an affine cone
over a projective orbifold.

\hfill

This is used to prove the following result
(for an alternative, more algebraic, proof, see 
\ref{_Stein_completion_second_proof_Proposition_}).

\hfill

\proposition\label{_Stein_completion_cone_Proposition_}
Let $P$ be a projective orbifold, $\dim_\C P >1$, and $L$ 
an ample line bundle on $P$. Assume that 
the total space $\Tot^\circ(L)$
of all non-zero vectors in $L$
is smooth. Then a weak Stein completion
of $\Tot^\circ(L)$ is obtained by adding a
point (called ``the origin'' or ``the apex'' elsewhere).

\hfill

\pstep
Let $P$ be a principal $\C^*$-bundle over $X$, with $\C^*$-action
denoted by $\rho_P(t):\; P \arrow P$. Consider the
dual $\C^*$-bundle $P^*$, which has the same underlying 
complex manifold but the group $\C^*$ acts on $P^*$
via $\rho_{P^*}(t) =  \rho_{P}(t^{-1})$. 
There is a natural duality between $P$ and $P^*$
mapping $P\times_X P^*$ to a trivial $\C^*$-bundle $\calo^*(M)$,
taking a pair $(p, p)$ of sections of $P$ and $P^*$, identified with $P$,
to the unit section of $\calo^*(X)$.

Let $\Tot^\circ(L)$ be the set of  non-zero
vectors in an ample line bundle $L$ over a projective
orbifold. Clearly, its dual bundle can be identified with
$\Tot^\circ(L^*)$, hence the spaces 
$\Tot^\circ(L)$ and $\Tot^\circ(L^*)$
are biholomorphic. 

\hfill
 
{\bf Step 2:}
Replacing $L$ by $L^*$, we can assume that
$L$ is anti-ample and it is equipped with a metric
with negative curvature. Then
the function $v\arrow |v|^2$ is strictly plurisubharmonic
on $\Tot^\circ(L)$, by \cite[(2.6)]{_Calabi:hk_}
or  \cite[(15.19)]{besse}, and defines
an LCK potential on the compact manifold
$M:=\frac{\Tot^\circ(L)}{\langle  \gamma\rangle}$,
where $\gamma(x) = \lambda x$, $|\lambda| >1$.
Applying 
\ref{potcon} to $\Tot^\circ(L)$, considered as 
a $\Z$-covering of $M$, we obtain 
that the Stein completion of $\Tot^\circ(L)$
is obtained by adding a point in the origin.
By \ref{_normalization_of_weak_Stein_comple_Remark_}, $\tilde M_c$ is the
normalization of any weak Stein completion $Z$ of $\tilde M$.
However, the normalization map $\nu:\; \tilde M_c\arrow Z$
is continous, and identity outside of the origin, because
$Z$ is smooth outside of the origin. Therefore, $\nu$
 is a diffeomorphism.
\endproof

\hfill

We give a more algebraic proof of 
\ref{_Stein_completion_cone_Proposition_},
independent from the LCK geometry, in Section 
\ref{_Remmert_Section_}. Note that the argument
in Section \ref{_Remmert_Section_} is stronger,
because it works even when $\dim_\C P =1$.

\hfill

\ref{potcon} was used to prove the following Kodaira
type embedding theorem for LCK manifolds with potential.
Recall that {\bf a linear Hopf manifold}
is a manifold $H:= \frac{\C^n \backslash 0}{\langle A
  \rangle}$,
where $A\in \GL(\C^n)$ is a linear contraction.
In \cite{ov_pams} it was shown that all linear Hopf
manifolds admit an LCK metric with potential.

\hfill

\theorem (\cite{ov_lckpot,ov_indam})
\label{_Embedding_Theorem_} 
Let $(M,I,\omega)$
be a compact LCK manifold with potential, 
$\dim_\C M\geq3$. Then $(M,I)$ admits a holomorphic embedding to a
linear Hopf manifold. \endproof

\hfill

\remark
The converse assertion is also clearly true:
a submanifold of a linear Hopf manifold is LCK with
potential (\cite{ov_pams}). 
Therefore, in dimension $\geq 3$, an LCK manifold with
(proper) potential can be defined as a manifold
which admits a holomorphic embedding into a linear
Hopf manifold. In dimension 2, this is also true,
if we assume the Kato conjecture (also known as
``the global spherical shell conjecture''); see
\cite[Section 5]{ov_indam} for details.


\section{The total space of an ample bundle, Remmert-Stein theorem 
and weak Stein completions}
\label{_Remmert_Section_}


In this section, we give an alternative
proof of \ref{_Stein_completion_cone_Proposition_},
independent from the LCK geometry. We start by
recalling the Remmert-Stein theorem.

\hfill

\theorem\label{_Remmert_Stein_Theorem_}
(Remmert-Stein theorem)\\ 
Let  $B$ be a complex analytic variety, and $C \subset B$ and 
 $A \subset B \backslash C$
 complex analytic subvarieties. Assume that
all irreducible components $A_i$ of $A$ satisfy
$\dim A_i > \dim C$. Then the closure of $A$ is
complex analytic in $B$.

\proof
\cite[\S II.8.2]{demailly}.
\endproof

\hfill

\remark \label{_Stein_Remmert_completion_Remark_}
We will use \ref{_Remmert_Stein_Theorem_}
in one special case, which is also the case used
in the proof of Chow theorem (\cite[Theorem II.8.10]{demailly}). Suppose that
$A \subset \C^n \backslash 0$ is a complex
analytic subset without 0-dimensional connected components; then
its closure $\bar A$ in $\C^n$ is complex analytic. A posteriori,
$\bar A$ is isomorphic to a weak Stein completion of
$A$, \ref{_Stein_completion_Definition_}.

\hfill

At this point we can 
give an alternative proof of \ref{_Stein_completion_cone_Proposition_}.

\hfill

\proposition\label{_Stein_completion_second_proof_Proposition_}
Let $P$ be a projective orbifold, and $L$ 
an ample line bundle on $P$. Assume that 
the total space $\Tot^\circ(L)$
of all non-zero vectors in $L$
is smooth. Then a weak Stein completion
of $\Tot^\circ(L)$ is obtained by adding a
point (called ``the origin'' or ``the apex'' elsewhere).
Moreover, the natural normalization
map from the Stein completion to a weak Stein 
completion is a homeomorphism.

\hfill

{\bf Proof:} The space
$\Tot^\circ(L)$ can be interpreted
algebraically as follows. 
Consider the standard $\C^*$-action
on $\Tot^\circ(L)$ understood
as a principal $\C^*$-bundle.
A finite-dimensional $\C^*$-representation
is a direct sum of 1-di\-men\-si\-o\-nal irreducible representations,
with $\C^*$ acting by $\rho(t) = t^w$; the number $w$ is called
{\bf the weight} of the representation. Therefore, the
category of $\C^*$-representations is equivalent
to the category of graded vector spaces. A 
ring with $\C^*$-action is the same as a graded ring.
In particular, the ring of fiberwise polynomial
holomorphic functions on $\Tot^\circ(L)$ is graded,
with the functions of weight $w$ being polynomials of
degree $w$ on each $\C^*$-orbit. We identify
the space $(\calo_{\Tot^\circ(L)})_w=(\calo_{\Tot^\circ(L^*)})_w$ of functions
of degree $w$ with the space of sections of $L^{\otimes w}$.
The variety $P\subset {\Bbb P} (H^0(L^N)^*)$ 
is identified with the graded spectrum
of the graded ring 
$\bigoplus_w H^0(X, L^{\otimes wN})$, whenever $L^N$ is very ample.
Then, the manifold $\Tot^\circ(L)$ is the preimage of $P$ under the
natural projection 
$H^0(L^N)^*\backslash 0 \arrow {\Bbb P} (H^0(L^N)^*)$.
The Stein completion of $\Tot^\circ(L)$
is its closure in $H^0(L^N)^*$, obtained by 
adding the zero (called ``the origin'' or 
``the apex'' elsewhere). It is
complex analytic by Remmert-Stein theorem
(\ref{_Stein_Remmert_completion_Remark_}), and normal by
\ref{proj_normality_Proposition_} below.

By \ref{_normalization_of_weak_Stein_comple_Remark_},
the Stein completion of $\Tot^\circ(L)$ is the
normalization of any weak Stein completion $Z$ of $\tilde M$.
Since the normalization map $\nu:\; \Tot^\circ(L)\arrow Z$
is continous, and identity outside of the origin, $\nu$
it is a homeomorphism.
\endproof


\section{Montel theorem and Riesz-Schauder theorem}


In this section, we recall several basic notions of 
complex analysis and functional analysis. 

\subsection{Montel theorem}
\definition \label{normal_family}
Let $M$ be a complex manifold, 
and ${\cal F}\subset H^0(\calo_M)$ a 
family of holomorphic functions. We call ${\cal F}$ 
{\bf a normal family} if for each compact
$K\subset M$ there exists $C_K>0$ such that
for each $f\in {\cal F}$, $\sup_K |f| \leq C_K$.

\hfill

\definition 
The {\bf $C^0$-topology} on the space
of functions on $M$ is the topology
of uniform convergence on compacts.

\hfill

\theorem \label{_Montel_Theorem_}
({\bf Montel})\\
Let ${\cal F}\subset H^0(\calo_M)$
be a normal family of functions,
and $\bar {\cal F}\subset H^0(\calo_M)$ its closure
in the $C^0$-topology. Then $\bar {\cal F}$ is compact
in $C^0$-topology.

\proof 
\cite[Lemma 1.4]{_Wu:Montel_}. \endproof

\subsection{The Banach space of holomorphic functions}

We briefly recall some notions of functional analysis, which are
standard.

\hfill

\remark
Let $C^0(X)$ be the space of continuous functions on a compact
space, with the norm $\|f\|_\sup:= \sup_{x\in X} |f(x)|$. It is not
hard to see that this space is Banach. Similarly,
if $X$ is a compact smooth manifold equipped with a connection $\nabla$,
the space $C^k(M)$ of $k$ times differentiable functions 
with the norm $\| f\|_{C^k}:=\sum_{i=0}^k \|\nabla^i(f)\|_\sup$
is also Banach.

\hfill

\theorem\label{_Banach_bounded_holo_Theorem_}
Let $M$ be a complex manifold, and let $H^0_b(\calo_M)$  
the space of all bounded holomorphic functions, equipped
with  the sup-norm $|f|_{\sup} := \sup_M |f|$.
 Then $H^0_b(\calo_M)$ is a Banach space.

\hfill

\proof 
Let $\{f_i\}\in H^0_b(\calo_M)$ be a Cauchy sequence in
the $\sup$-norm. {Then $\{f_i\}$ converges to a continuous
	function $f$} in the $\sup$-topology.

Since $\{f_i\}$ is a normal family (see \ref{normal_family}), it
has a subsequence which converges in $ C^0$-topology
to $\tilde f\in H^0(\calo_M)$, by Montel's \ref{_Montel_Theorem_}. 
However, the $ C^0$-topology
	is weaker than the $\sup$-topology, hence $\tilde f=f$.
Therefore, $f$ is holomorphic. \endproof

\subsubsection{Compact operators}

Recall that 
a  subset $X$ of a topological space $Y$ is called {\bf 
	precompact},
or {\bf relatively compact in $Y$}, 
if its closure is compact. 

\hfill

\definition \label{_bounded_set_Definition_}
A subset
$K\subset V$ of a topological vector space
is called {\bf bounded}   if for any 
open set $U\ni 0$, there exists a number $\lambda_U\in
\R^{>0}$ such that $\lambda_U K \subset U$.

\hfill

\definition \label{_compact_operator_Definition_}
Let $V, W$ be topological vector spaces. A continuous operator $\phi:\; V
\arrow W$ is called {\bf  compact} if the image of
any bounded set is precompact.

\subsubsection{Holomorphic contractions}

\theorem \label{_contra_compact_Theorem_}
Let $X$ be a complex variety, and 
$\gamma:\; X \arrow X$ a holomorphic contraction
to $x\in X$ such that $\gamma(X)$ is precompact. 
Consider the Banach space $V=H^0_b(\calo_X)$ of bounded holomorphic
functions with the sup-norm, and let $V_x\subset V$
be the space of all $v\in V$ vanishing in $x$. Then $\gamma^*:\; V \arrow V$
	is compact, and the eigenvalues of its restriction to 
	$V_x$ are strictly less than 1 in absolute value.\footnote{Since
$\gamma^*$ maps constants to constants identically, one cannot expect
that $\|\gamma^*\|< 1$ on $V$. However, if we add a condition
which excludes constants, such as $v(x)=0$, we immediately
obtain $\|\gamma^*\|< 1$.}

\hfill

\pstep For any $f\in H^0(\calo_X)$ we have
\[|\gamma^* f|_{\sup}= \sup_{x\in \overline{\gamma(X)}}
|f(x)|.
\]
This implies that $\gamma^*(f)$ is bounded.
Therefore, {for any sequence $\{f_i\in H^0(\calo_X)\}$ converging in the 
	$ C^0$-topology, the sequence $\{\gamma^* f_i\}$ converges
	in the $\sup$-to\-po\-lo\-gy.}
The set $B_C:=\{v\in V\ \ |
\ \  |v|_{\sup} \leq C\}$ is precompact in the 
$ C^0$-topology, because it is a normal family.
Then $\gamma^* B_C$ is precompact in the 
$\sup$-topology.\footnote{The $C^0$-convergence for holomorphic
functions is strictly weaker that the $\sup$-convergence; for example, the sequence of bounded holomorphic
functions on an open disk, $f_i:= z^i$ converges in $C^0$-topology, but does not converge in the $\sup$-topology.} 
This proves that the operator $\gamma^*:\; V\arrow V$ is compact.

It remains to show that its operator norm is $<1$ on $V_x$.

\hfill

{\bf Step 2:}
Since $\sup_X |\gamma^* f|= \sup_{\gamma(X)} |f| \leq \sup_X |f|$,
one has $\|\gamma^*\|\leq 1$. If this inequality is not
strict, for some sequence $\{f_i\}$ of holomorphic
functions $f_i\in V_x$ with $\sup_X |f_i|\leq 1$ (that is, $f_1\in B_1$)
one has $\lim_i \sup_{x\in \gamma(X)} |f_i(x)|=1$.
Since $B_1$ is a normal family, $f_i$ has a subsequence
converging in $ C^0$-topology to $f$. Then $\{\gamma(f_i)\}$
converges to $\gamma(f)$ in $\sup$-topology, giving 
$$\lim_i \sup_{x\in \gamma(X)} |f_i(x)|= \sup_{x\in  \gamma(X)} |f(x)|=1.$$
Since $f(x)=0$, $f$ is non-constant.
By the maximum modulus principle,
a non-constant holomorphic function has no local maxima; this means
	that $|f(x)| >1$ somewhere on $X$. Then $f$ cannot
be the $C^0$-limit of $\{f_i\}\subset B_1$. \endproof

\subsubsection{The Riesz-Schauder theorem}

The following result  is a Banach analogue of 
the usual spectral theorem for compact operators on
Hilbert spaces. 

\hfill

\theorem ({Riesz-Schauder}, \cite[Section 5.2]{friedman})\\
\label{_Riesz_Schauder_main_Theorem_}
Let $F:\; V \arrow V$ be a compact operator on a Banach space.
Then for each non-zero $\mu \in \C$, there exists a sufficiently
big number $N\in \N$ such that for each $n>N$ {one has 
	$V= \ker(F-\mu\Id)^n \oplus \overline{\im (F-\mu\Id)^n}$,
	where $\overline{\im (F-\mu\Id)^n}$ is the closure of the image.}
Moreover, $\ker(F-\mu\Id)^n$ is finite-dimensional and independent on $n$.
\endproof

\hfill

\remark\label{_root_space_RS_Remark_}
Define  {\bf the root space of an operator 
$F\in \End(V)$, associated with an eigenvalue $\mu$,}
as $\bigcup_{n\in \Z} \ker(F-\mu\Id)^n$. 
In the finite-dimensional case, the root spaces
coincide with the  Jordan cells of the corresponding 
matrix. Then 
\ref{_Riesz_Schauder_main_Theorem_} can be reformulated
by saying that any compact operator $F\in \End(V)$ admits a 
Jordan cell decomposition, with $V$ identified with
a completed direct sum of the root spaces, which are
all finite-dimensional; moreover, the eigenvalues
of $F$ converge to zero.

\hfill

We need the following corollary of the Riesz-Schauder theorem, which is
obtained using the same arguments as in the finite-dimensional case.

\hfill

\theorem\label{rs}
Let  $F:\; V \arrow V$ be a compact operator on a Banach space.
We say that $v$ is {\bf a root vector} for $F$ if $v$ lies in 
a root space of $F$, for some eigenvalue $\mu\in \C$. 
Then the space generated by the root vectors is dense in $V$.
\endproof

\hfill

\definition\label{_F_finite_Definition_}
Let $F\in \End(V)$ be an endomorphism of a vector space.
A vector $v\in V$ is called {\bf $F$-finite}
if the space generated by $v, F(v), F(F(v)),  ...$
is finite-dimensional. 

\hfill

\remark Clearly, a vector is 
$F$-finite if and only if it is obtained as a combination
of root vectors. Then  \ref{rs} 
implies the following.

\hfill

\corollary\label{_finite_RS_Corollary_}
 Let  $F:\; V \arrow V$ be a compact operator on a Banach space,
and $V_0\subset V$ the space of all $F$-finite vectors.
Then $V_0$ is dense in $V$.
\endproof

\hfill

We also need the following lemma, giving ``a relative case''
of Riesz-Schauder.

\hfill

\lemma\label{_RS_surj_Lemma_}
Consider the commutative square of continuous 
operators of Banach spaces
\[ \begin{CD}
W_1 @> K_1 >> W_1\\
@V R VV @VV R V \\
W_2 @> K_2 >> W_2
\end{CD}
\]
with $R$ a surjective map, and both $K_i$, $i=1,2$, are compact.
Restricting $R$ to the space $W_i^{K_i}$ of $K_i$-finite vectors, we
obtain a map $R_f:\; W_1^{K_1}\arrow  W_2^{K_2}$.
Then $R_f$ is surjective.

\hfill

\proof
Let $U_2\subset W_2$ be a finite-dimensional space
of $K_2$-finite vectors, and $U_1\subset W_1$ its preimage.
Since $R$ is surjective, the natural map
$U_1 \to U_2$ is surjective. Since $U_1^{K_1}$ is dense in $U_1$,
and $U_2$ is finite-dimensional, the restriction
of $R$ to $U_1^{K_1}$ is also surjective.
\endproof


\section{Algebraic structures on weak Stein completions}
\label{_algebraic_comple_Section_}

In this section,
we put an algebraic structure
on a weak Stein completion of a K\"ahler $\Z$-covering
of a submanifold in a Hopf manifold.
We also prove that this algebraic structure is unique.
Later on, this algebraic structure is used to
show that this weak Stein completion is an algebraic
cone.

\subsection{Ideals of the embedding to a Hopf manifold}

Let $M\subset H$ be a subvariety of a linear Hopf manifold
$H=\frac{\C^n \backslash 0}{\langle A \rangle}$.
In this subsection we prove that 
the preimage of $M$ in $\C^n \backslash 0$ can be defined by a set 
of polynomial equations.

\hfill

We start with the following lemma.

\hfill

\lemma\label{_gamma_finite_on_C^n_Lemma_}
Let $\gamma$ be an invertible linear contraction of $\C^n$.
A holomorphic function on $\C^n$
is $\gamma$-finite if and only if it is polynomial.

\hfill

\proof
Clearly, a polynomial function is $\gamma$-finite.
The operator $\gamma^*$
acts on homogeneous polynomials of degree $d$
with eigenvalues $\alpha_{i_1}\alpha_{i_2}... \alpha_{i_d}$,
where $\alpha_{i_j}$ are the eigenvalues of $\gamma$ on $\C^n$.
Since $\gamma$ is a contraction, all $\alpha_{i_j}$
satisfy $|\alpha_{i_j}| <1$. Therefore,
any sequence $\{\alpha_{i_1}\alpha_{i_2}... \alpha_{i_d}\}$
converges to 0 as $d$ goes to infinity. 
We obtain that every given number can be
realized as an eigenvalue of $\gamma^*$ on
homogeneous polynomials of degree $d$ 
for finitely many choices of $d$ only.
Therefore, any root vector of $\gamma^*$
is a finite sum of homogeneous polynomials.
This implies that
the Taylor decomposition of a $\gamma$-finite function $f$ can   have 
only finitely many components, otherwise the eigenspace
decomposition of $f$ with respect to the
action of $\gamma^*$ is infinite. 
\endproof

\hfill

This lemma can be used to prove the following theorem
(see also \ref{_algebraic_str_on_cone_existence_Theorem_} below).

\hfill

\theorem \label{_cone_cover_for_LCK_pot_algebraic_Theorem_}
{ (\cite[Theorem 2.8]{ov_pams}, \cite[Theorem 5.5]{ov_lee}})\\
Let $M \hookrightarrow H$ be a complex subvariety of a linear Hopf manifold,
and $\tilde M_c \arrow \C^N$ the corresponding map of weak Stein completions,
with $\tilde M_c$ obtained as the closure of $\tilde M \subset \C^N$
by adding the zero
(\ref{_Stein_Remmert_completion_Remark_}). 
Then $\tilde M_c$ is an algebraic subvariety, that is,
a set of common zeroes of a system of polynomial equations.
\endproof

\subsection{Algebraic structures on Stein completions: the existence}

Let $X\subset M \subset \C P^n$ be
projective subvarieties of $\C P^n$.
Then $Z= M\backslash X$ is called {\bf a quasi-projective variety}.
The {\bf Zariski topology} on $Z$ is the topology such that the
closed subsets of $Z$ are closed algebraic subvarieties.
The {\bf regular functions on a Zariski open subset
  $W\subset X$} 
are restrictions of the
rational functions on $\C P^n$ which have no
poles on $W$. This defines 
 {\bf the sheaf of regular functions} on $Z$
as a subsheaf of the sheaf of holomorphic functions
in the Zariski topology on $Z$.

Recall that {\bf an algebraic structure} on a 
complex analytic variety $Z$ is a subsheaf of the
sheaf of holomorphic functions which can be
realized as a sheaf of regular functions 
for some biholomorphism between $Z$ and
a  quasi-projective variety.

Note that the algebraic structure is not
unique, indeed, even the complex manifold
$(\C^*)^2$ has more than one algebraic structure
(\cite{_Simpson:rank_one_}). In \cite{_Jelonek:exotic_},
Z. Jelonek has constructed an uncountable set 
of pairwise non-isomorphic algebraic structures on
certain Stein manifolds.

It turns out that 
a K\"ahler $\Z$-cover of an LCK manifold with potential
is equipped with a distinguished (affine) algebraic structure, which is
uniquely determined by the complex structure
(\ref{_algebraic_str_on_cone_existence_Theorem_},
\ref{_same_algebra_structure_on_cone_Theorem_}).

\hfill

The main result of this section is the following theorem.

\hfill

\theorem\label{_algebraic_str_on_cone_existence_Theorem_}
Let $M$ be an LCK manifold with proper potential,
and  $\tilde M_c$ the Stein variety
obtained as a weak Stein completion of its K\"ahler $\Z$-cover $\tilde M$.
Denote by $\gamma\in \Aut(\tilde M_c)$ 
the generator of the $\Z$-action on $\tilde M$ extended to 
$\tilde M_c$. This map is an invertible holomorphic contraction
centered in the apex $c$ of $\tilde M_c$. 
Then $\tilde M_c$ admits an algebraic structure such that
the regular functions are precisely the $\gamma$-finite functions.

\hfill

\pstep
\ref{_gamma_finite_on_C^n_Lemma_}
implies \ref{_algebraic_str_on_cone_existence_Theorem_}
when $\gamma$ is a linear contraction and $\tilde M_c=\C^n$.

\hfill

{\bf Step 2:} By definition, $M$ admits a holomorphic
embedding to a linear Hopf manifold. 
Consider the corresponding embedding 
$\tilde M_c \hookrightarrow \C^n$.
By \ref{_cone_cover_for_LCK_pot_algebraic_Theorem_}, 
$\tilde M_c$ is an algebraic subvariety of $\C^n$ 
given by a finite set of polynomial equations.

To finish the proof of \ref{_algebraic_str_on_cone_existence_Theorem_},
it remains to show that a holomorphic function on $\tilde M_c$ is 
polynomial if and only if it is $\gamma$-finite.
It is clear that any holomorphic function on $\tilde M_c$
obtained from a polynomial function on $\C^n$
is $\gamma$-finite. It remains to show the converse.
Consider a $\gamma$-finite function on $\tilde M_c$. By 
\ref{_gamma_finite_on_C^n_Lemma_},
to prove that it is polynomial, it would suffice
to prove that it can be extended to 
a $\gamma$-finite function on $\C^n$.

Consider the exact sequence of sheaves of holomorphic
functions on $\C^n$
\[ 
0 \arrow I_{\tilde M_c} \arrow \calo_{\C^n}\arrow \calo_{\tilde M_c}\arrow 0,
\]
where $I_{\tilde M_c}$ is the ideal sheaf of $\tilde M_c$.
The corresponding long exact sequence gives
\[
0 \arrow H^0(\C^n,I_{\tilde M_c}) \arrow H^0(\C^n,\calo_{\C^n})\arrow 
H^0(\C^n,\calo_{\tilde M_c})\arrow H^1(\C^n,I_{\tilde M_c}).
\]
Since $\C^n$ is Stein, $H^1(\C^n,I_{\tilde M_c})=0$,
hence the function $f$ is the restriction of
a holomorphic function $\tilde f \in H^0(\C^n,\calo_{\C^n})$.
Finally, by \ref{_RS_surj_Lemma_}, $\tilde f$ can be chosen
$\gamma$-finite.
\endproof

\subsection{Algebraic structures on Stein completions: the uniqueness}

In this subsection, we prove that the 
algebraic structure constructed on $\tilde M_c$ in 
\ref{_algebraic_str_on_cone_existence_Theorem_}
is uniquely determined by the complex geometry of the
Stein variety $\tilde M_c$.

We start with the following lemma, showing that
a function is $\gamma$-finite (that is, polynomial)
if and only if it admits a certain growth condition.
This is similar to a well-known result of complex
analysis which states that an entire holomorphic
function on $\C^n$ is polynomial if and only if
it has polynomial growth. 

\hfill

\lemma\label{_poly_growth=gamma_finite_Lemma_}
Let $M$ be an LCK manifold with proper potential,
and  $\tilde M_c$ the Stein variety,
obtained as a weak Stein completion of its K\"ahler 
$\Z$-cover $\tilde M$. Denote by $\gamma\in \Aut(\tilde M_c)$ the
holomorphic contraction generating the $\Z$-action on
$\tilde M_c$. Choose a compact set $K\subset \tilde M_c$
containing an open neighbourhood of the apex. Denote by
${\cal B}\subset H^0(\calo_{\tilde M_c})$ 
the following ring of functions on $\tilde M_c$:
\begin{equation}\label{_polynomial_growth_Equation_}
{\cal B}:= \{ f \in H^0(\calo_{\tilde M_c})\ |\ 
\exists\; C>0 \text{\ such that\ } \forall i \ \ \sup_{K}|(\gamma^*)^{-i}|f| < C^i\}.
\end{equation}
Then ${\cal B}$ coincides with the space of $\gamma$-finite 
functions.\footnote{We call the function satisfying
\eqref{_polynomial_growth_Equation_} {\bf the function of polynomial growth}.
This terminology is justified because for $\gamma$ a linear contraction 
of $\C^n$, \eqref{_polynomial_growth_Equation_} is equivalent to having 
polynomial growth.}

\hfill

\pstep
Let $f$ be a $\gamma$-finite function, and $W$ the
space generated by $\{f, (\gamma^*)f, (\gamma^*)^2 f, ...\}$.
Let $\|\cdot \|_K$ be the norm on $W$ defined by
$\|f\|:= \sup_{K} |f|$, and let
$C:= \sup_{\|f\|=1}\|(\gamma^*)^{-1}f\|$ be
the operator norm of the operator
$(\gamma^*)^{-1}\in \End(W)$ in this norm.
Then $\sup_{K}|(\gamma^*)^{-i}f|\leq C^i \sup_{K} |f|$,
and $f$ has polynomial growth, in the sense of 
\eqref{_polynomial_growth_Equation_}.

\hfill

{\bf Step 2:} 
Suppose that $f$ has polynomial growth, in the sense 
of \eqref{_polynomial_growth_Equation_},
and let $W$ be the space of all functions
generated by $\{f, (\gamma^*)f, (\gamma^*)^2 f, ...\}$.
Then all elements of $W$ have the same growth as $f$,
with the same bound $C$, hence the closure 
$\bar W$ of $W$ in the norm $\|f\|:= \sup_{K} |f|$ consists
of functions with polynomial growth. 

The condition \eqref{_polynomial_growth_Equation_}
holds for all $f \in W$ if and only if 
$(\gamma^*)^{-1}$ has finite norm on
$W$. Therefore, $\gamma^*\restrict{\bar W}$ is invertible. 
To finish the proof of \ref{_poly_growth=gamma_finite_Lemma_},
we need to prove that the norm of $(\gamma^*)^{-1}$
is infinite on $\bar W$ if $\bar W$ is infinite-dimensional.

The operator $\gamma^*$ on $\bar W$ is compact;
by the Riesz-Schauder theorem, it has the Jordan cell
decomposition with eigenvalues converging to 0, unless
$\bar W$ is finite-dimensional.
The norm of a linear operator $A$ with eigenvalues
$\alpha_i$ satisfies $\|A\| \geq \sup |\alpha_i|$.
Therefore, a compact operator cannot be 
invertible on an infinitely-dimensional Banach space:
the inverse operator would have infinite norm.
\endproof

\hfill

\theorem\label{_same_algebra_structure_on_cone_Theorem_}
Let $\tilde M_c$ be a weak Stein completion of a K\"ahler
$\Z$-covering of an LCK manifold with potential.
Suppose that $\tilde M_c$ can be obtained 
from two different LCK manifolds, and let
$\gamma_1, \gamma_2 \in \Aut(\tilde M_c)$ be
the corresponding holomorphic contractions. 
Consider the affine algebraic structures ${\cal A}_1$ and ${\cal A}_2$
associated with the LCK structures
on $\frac{\tilde M}{\langle \gamma_1\rangle}$
and $\frac{\tilde M}{\langle \gamma_2\rangle}$
as in \ref{_algebraic_str_on_cone_existence_Theorem_}.
Then ${\cal A}_1={\cal A}_2$.

\hfill

\pstep
Let $\gamma_1, \gamma_2 \in \Aut(\tilde M_c)$
be holomorphic contractions centered in the
apex $c$, associated with two LCK manifolds
with proper potential such that their  K\"ahler 
$\Z$-covers have the same weak Stein
completion $\tilde M_c$.
We prove \ref{_same_algebra_structure_on_cone_Theorem_}
using \ref{_poly_growth=gamma_finite_Lemma_}.
To show that the algebraic structures induced by
$\gamma_1$ and $\gamma_2$ are the same, it would suffice
to show that the spaces of $\gamma_1$- and $\gamma_2$-finite
functions on $\tilde M_c$ coincide. 
By \ref{_poly_growth=gamma_finite_Lemma_}, 
this would follow
if the growth estimates 
\eqref{_polynomial_growth_Equation_}
for $\gamma_1$ and $\gamma_2$ are equivalent.

We have reduced \ref{_same_algebra_structure_on_cone_Theorem_}
to the following statement. Let $K\subset \tilde M_c$
be a compact subset which contains an open neighbourhood of 
the apex $c\in \tilde M_c$. Then 
a function $f \in H^0(\calo_{\tilde M_c})$ which satisfies
$\sup_K(\gamma_1)^{-n}|f| < C_1^n$, also
satisfies  $\sup_K(\gamma_2)^{-n}|f| < C_2^n$
for some other constant $C_2 >0$ and all $n>0$.

\hfill

{\bf Step 2:} Suppose that there
exists an integer $d>0$ such that
$\gamma_1^d(K) \subset \gamma_2(K)$.
Then the polynomial growth estimate \eqref{_polynomial_growth_Equation_}
for $\gamma_2$ follows from the growth estimate for
$\gamma_1^d$, which is equivalent to the 
growth estimate for $\gamma_1$.
Therefore, \ref{_same_algebra_structure_on_cone_Theorem_}
would follow if we prove that the integer $d$
always exists.

Recall that a continuous map $f:\; M \arrow M$ fixing $x\in M$
is called {\bf a contraction centered in $x\in M$}
if for each compact subset $K\subset M$, 
and each open neighbourhood of $x$, a sufficiently
big iteration of $f$ gives $f^N(K) \subset U$.
The maps $\gamma_1, \gamma_2$ are contractions centered in 
the apex $x$, hence for some $N>0$, the set 
$\gamma_1^N(K)$ lies in the interior
of $\gamma_2(K)$ which is a neighbourhood of $x$.
\ref{_same_algebra_structure_on_cone_Theorem_} is proven.
\endproof


\section{Algebraic cones and LCK manifolds}
\label{_cones_and_LCK_Section_}

In Section \ref{_cones_Section_},
we gave three equivalent definitions
of an algebraic cone: an open algebraic
cone is the total space of an ample $\C^*$-bundle,
or a K\"ahler $\Z$-covering of a submanifold
in a Hopf manifold, or a Stein variety admitting a
holomorphic contraction, with the apex removed.
In this section, we prove this equivalence.
We start by recalling some basic facts 
about the weighted projective spaces.

\subsection{Weighted projective spaces}

For an introduction and background
material about the weighted projective spaces,
see \cite[\S 4.5]{bog}. 

\hfill

Recall that any representation $V$ of $\C^*$
is a direct sum of 1-dimensional representations isomorphic to $\rho_w$,
with $\C^*$ acting by $\rho_w(t)(z)= t^w z$.
Such a representation is called {\bf representation of weight $w$},
and the decomposition of $V$ into subrepresentation 
of constant weight -- {\bf the weight decomposition}.

\hfill

The following trivial 
claim is left as an exercise to the reader.

\hfill

\claim\label{_all_weights_positive_Claim_}
Let $\rho$ be $\C^*$ acting on $\C^n$.
Assume that $\rho$ contains a contraction.
Then all weights of $\rho$ are positive or negative.
\endproof

\hfill

\claim\label{_weighted_proj_Claim_}
Let $\rho$ be $\C^*$ acting on $\C^n$
with weights $w_1, ..., w_n\in \Z^{>0}$.
Then its orbit space $\C P^{n-1}(w_1, ..., w_n)$ is 
equipped with a structure of a projective
orbifold, and $\C^n \backslash 0$
can be identified with the total space
of an ample $\C^*$-bundle over $\C P^{n-1}(w_1, ..., w_n)$.

\proof \cite[Proposition 4.5.3]{bog}.
\endproof

\hfill

\definition
The orbifold $\C P^{n-1}(w_1, ..., w_n)$
is called {\bf the weighted projective space}.

\hfill

\claim\label{_subvariety_in_weighted_Claim_}
Let $\rho$ be $\C^*$ acting on $\C^n$
with weights $w_1, ..., w_n\in \Z^{>0}$,
and $Z \subset \C^n\backslash 0$ be a $\rho$-invariant submanifold.
Then the orbit space $Z/\C^*$ is a projective orbifold
in the corresponding weighted projective space
$\C P^{n-1}(w_1, ..., w_n)$.

\proof \cite[Proposition 4.6.2]{bog}.
\endproof

\hfill

\corollary\label{_C^*_action_then_cone_Corollary_}
Let $\rho$ be $\C^*$ acting on $\C^n$
with weights $w_1, ..., w_n\in \Z^{>0}$,
and  $Z \subset \C^n\backslash 0$ a $\rho$-invariant submanifold.
Then $Z$ is an open algebraic cone.
\endproof

\subsection{Jordan-Chevalley decomposition}

We recall a few definitions and results from the algebraic group 
theory, generalizing the Jordan normal form to arbitrary
algebraic group.

\hfill

\definition 
An element of an algebraic group $G$ is called
{\bf semisimple} if its image is semisimple
for some exact
algebraic representation of $G$, and is called 
{\bf unipotent} if its image is unipotent
(that is, exponential of a nilpotent) 
for some exact algebraic representation of $G$.

\hfill

\remark 
For any algebraic representation of an algebraic group
$G$, the image of any semisimple element is a 
semisimple operator, and the image of any unipotent
element is a unipotent operator (\cite[\S 15.3]{hum}).

\hfill

\theorem  ({Jordan-Chevalley  decomposition}, 
\cite[\S 15.3]{hum}) \label{jcdec}\\
Let $G$ be an algebraic group, and $A\in G$.
{Then there exists a unique decomposition $A= S U$
	of $A$ in a product of commuting elements $S$ and $U$,
	where $U$ is unipotent and $S$ semisimple.}
\endproof

\subsection{Algebraic cone obtained from an LCK manifold with
potential}

\theorem\label{_cone_cover_for_LCK_pot_Theorem_}
Let $M$ be an LCK manifold with proper potential,
and $\tilde M$ its K\"ahler $\Z$-cover.
Then $\tilde M$ is an open algebraic cone.

\hfill

\proof
Let $\gamma$ be the generator of the $\Z$-action on 
$\tilde M$.
Using \ref{_Embedding_Theorem_}, we can 
embed $M$ to a linear Hopf manifold.
This gives a holomorphic map
$\tilde M \arrow \C^n\backslash 0$ taking
$\gamma$ to a linear contraction $A\in \End(\C^n)$.
Let $\tilde M_c$ be the closure of $\tilde M$
in $\C^n$. By \ref{_cone_cover_for_LCK_pot_algebraic_Theorem_},
$\tilde M_c$ is equipped with a natural algebraic structure.
To prove \ref{_cone_cover_for_LCK_pot_Theorem_}
it remains only to show that 
this variety admits a holomorphic $\C^*$-action, free
outside of the origin $c$, and acting
with eigenvalues $|\alpha_i|< 1$ on 
the Zariski tangent space $T_c \tilde M_c$.

Let ${\cal G}_A$ be the Zariski closure of
$\langle A \rangle$ in $\GL(\C^n)$. This is a commutative
algebraic group, acting on the variety $\tilde M_c\subset \C^n$.
To prove that $\tilde M_c$ is an algebraic cone, we find
a $\C^*$-action contracting $\tilde M_c$ to the origin and
apply \ref{_C^*_action_then_cone_Corollary_}.

Let $A=SU$ be the Jordan-Chevalley decomposition for
$A$, with $S, U\in {\cal G}_A$. Since ${\cal G}_A$ preserves $\tilde M_c$, 
the endomorphisms $S$ and $U\in \End(\C^n)$ also act on $\tilde M_c$.
Since the eigenvalues of $S$ are the same
as the eigenvalues of $A$, the map $S$ is a linear contraction.
Let $G_S\subset {\cal G}_A$, $G_S= e^{\C \log S}$ be a one-parametric subgroup
containing $S$. We prove that $G_S$ can be approximated by
subgroups of ${\cal G}_A$ isomorphic to $\C^*$; then these subgroups
also contain a contraction, and we can apply 
\ref{_C^*_action_then_cone_Corollary_}.

Consider the map taking any $A_1\in {\cal G}_A$ to its unipotent component
$U_1$. Since ${\cal G}_A$ is commutative, this map is a group
homomorphism. Therefore, its kernel ${\cal G}_s$ (that is,
the set of all semisimple elements in ${\cal G}_A$) is an
algebraic subgroup of ${\cal G}_A$. A connected commutative
algebraic subgroup of $\GL(\C^N)$ consisting of semisimple elements
is always isomorphic to $(\C^*)^k$ 
(\cite[Proposition 1.5]{_Borel_Tits:Groupes_Reductifs_}).
 The one-parametric subgroups $\C^*\subset (\C^*)^k$ 
are dense in $(\C^*)^k$ because one-parametric
complex subgroups $\C^*\subset (\C^*)^k$ can be obtained as 
complexification of subgroups $S^1\subset \U(1)^k \subset (\C^*)^k$,
and those are dense in $\U(1)^k$. Therefore,
the contraction $S\in {\cal G}_s=(\C^*)^k$ can be approximated
by an element of $\C^*$ acting on $\tilde M_c$.

We have obtained a $\C^*$-action $\rho$ on $\tilde M_c$
containing a contraction. By  \ref{_all_weights_positive_Claim_},
all weights of this $\C^*$-action are positive (or negative,
in which case we replace $\rho$ by its opposite). Then,
\ref{_C^*_action_then_cone_Corollary_} implies that $\tilde M$
is an open algebraic cone.
\endproof

\subsection{Submanifolds of Hopf manifolds obtained from
  algebraic cones}

\lemma\label{_properly_discont_contra_Lemma_}
Let $\gamma$ act on a complex variety 
$\tilde M_c$ by contractions, contracting $\tilde M_c$ to a point $c$.
Then the corresponding $\Z$-action on $\tilde M:= \tilde M_c\backslash \{c\}$ 
is properly discontinuous, hence $\tilde M/\Z$ is Hausdorff; it is 
a manifold when $\tilde M$ is a manifold.

\hfill

\proof
By definition, the $\Z$-action is properly discontinuous if every
point has a neighbourhood $U$ such that the set 
$\{g\in \Z \ \ |\  \ g(U)\cap  U \neq \emptyset\}$
is finite. Let $x\in \tilde M$ and 
$K$ be the compact closure of an open neighbourhood
of $U\subset \tilde M$ containing $x$. Since $\tilde M_c$ is Hausdorff,
there exists a neighbourhood $W\ni c$ such that
its closure does not intersect $K$. By definition of contractions,
there exists $N>0$ such that 
$\gamma^n(K)\subset W$ for all $n\geq N$.
This implies that $\gamma^n(K)\cap K = \emptyset$.
This also implies that $K \cap \gamma^{-n}(K)=\emptyset$.
We have shown that $\gamma^n(K)\cap K= \emptyset$
for all $n \notin [-N, N]$.
\endproof

\hfill

The following theorem, applied to a weak Stein
completion of an open algebraic cone, allows
us to obtain a holomorphic embedding of its
$\Z$-quotient to a Hopf manifold.

\hfill

\theorem\label{_alge_cone_to_Hopf_Theorem_}
Let $\tilde M_c$ be a Stein variety equipped with
a holomorphic contraction $\gamma$, contracting
it to the point $c$. Assume that the complement
$\tilde M:= \tilde M_c\backslash c$ is smooth.
  By \ref{_properly_discont_contra_Lemma_},
$\tilde M/\langle \gamma \rangle$ is a complex manifold.
Then there exists a holomorphic embedding 
$j:\; \tilde M/\langle \gamma \rangle\hookrightarrow H$
to a Hopf manifold.  

\hfill

\proof
Let $R$ be the ring of $\gamma$-finite holomorphic functions on
$\tilde M_c$, $I$ the maximal ideal of $c$, and
and $V\subset I$ a finite-dimensional 
$\gamma$-invariant space generating $I$. 
By \ref{_contra_compact_Theorem_}, 
the action of $\gamma^*$ is compact
on $I$ and has all eigenvalues $< 1$.
By Riesz-Schauder theorem, 
$R$ is dense in $\calo_{\tilde M_c}$ 
(\ref{_finite_RS_Corollary_}). Therefore, the
functions in $V$ separate the points in $\tilde M$,
for $V$ sufficiently big. 

By Montel's theorem, $R$ is dense in $C^1$-topology whenever it is
dense in $C^0$-topology; therefore, 
the differentials of the functions $f\in R$ generate $T^*_x\tilde M$
for all $x\in \tilde M$. This implies that
the tautological map $\tilde M \arrow V^*$
is a holomorphic embedding.
This map is by construction $\gamma$-equivariant,
hence it induces a holomorphic embedding
$\tilde M/\langle \gamma \rangle\hookrightarrow H=V^*/\langle \gamma\rangle$.
\endproof

%
%
%
%
%
%
%
%
%
%


\section{Closed algebraic cones and normal varieties}
\label{_normality_of_cone_Subsection_}

\definition\label{_normality_Definition_}
Recall that a complex variety $X$ is called {\bf normal}
if any locally bounded meromorphic function on an open subset 
$U\subset X$ is holomorphic (\cite[Definition II.7.4]{demailly}). 
In algebraic geometry, a variety
is normal if all its local rings are integrally closed; these
two notions are equivalent for a complex variety obtained
from an algebraic one (\cite[Theorem II.7.3]{demailly},
\cite[Satz 4, p. 122]{_Kuhlmann_}).\footnote{We are grateful to Francesco Polizzi
for this reference; please see the 
excellent Mathoverflow thread 
{\tiny\url{ https://mathoverflow.net/questions/303406/algebraic-vs-analytic-normality/}}
for more details about the complex analytic and complex algebraic normality.}

\hfill

\remark
Each complex variety $A$ admits {\bf the normalization},
that is, a normal variety $B$ equipped with a finite
holomorphic map $n:\; B \arrow A$. The normalization
is unique up to an isomorphism. 
In the sequel, we will usually consider 
the normalization maps $n:\; B \arrow A$
which are homeomorphisms. The concept might
seem alien, but in fact it is very natural.
For example, consider the plane curve $C$
defined by the equation $x^2=y^3$. This
curve is clearly singular in zero.
It can be parametrized by a variable
$t$ with $t^2=y$ and $t^3=x$; this
parametrization defines a map $n:\; \C \arrow C$
which is finite, hence $n$ is the
normalization of $C$. However, $n$ is clearly 
a homeomorphism.

\hfill

\claim\label{_quotient_normal_Claim_}
Let $G$ be a finite group holomorphically acting 
on a normal complex variety $Y$, and $X= Y/G$ the quotient 
variety. Then $X$ is normal.

\hfill

\proof Let $f$ be a locally bounded meromorphic function
on $X$. Its pullback $\tau^* f$ to $Y$ is clearly meromorphic and locally
bounded, hence holomorphic. Since $\tau^* f$ is $G$-invariant,
it is holomorphic on $Y$.
\endproof

\hfill

We are interested in the normality of a closed algebraic cone.

Let now $X$ be a projective variety and $L$ an ample line bundle
on $X$. The {\bf homogeneous coordinate ring} of $X$ is
the ring $\bigoplus_i H^0(X, L^{\otimes i})$. This ring is finitely
generated whenever $L$ is ample (\cite[Theorem 2.3.15]{_Lazarsfeld:1_}). 
Indeed, for an ample line bundle
$L$ and $i$ sufficiently big, the cohomology $H^{>0}(X,  L^{\otimes i})$ vanishes,
hence the number $\dim H^0(X,  L^{\otimes i})$ is equal to the holomorphic Euler 
characteristic $\chi(L^{\otimes i})$. The latter is polynomial in $i$
by Riemann-Roch-Hirzebruch theorem, which implies that finitely
many generators are enough to generate this ring.

A projective variety $(X, L)$ for which the ring
$\bigoplus_i H^0(X, L^{\otimes i})$ is integrally closed
is called {\bf projectively normal}. By definition, the
affine variety associated with this ring
is {\bf the affine cone of $X$}, which is a closed algebraic cone in our 
parlance. 

The closed algebraic cone, obtained from an open algebraic
cone by taking the Stein completion, is normal by
definition of Stein completion.
However, the projective normality is a very tricky condition,
depending on the choice of the bundle $L$.

This was the subject of Exercise I.3.18 
from Hartshorne's ``Algebraic Geometry'' (\cite{_Hartshorne:AG_}),
where Hartshorne considers the same smooth complex curve with two different
projective embeddings; the affine cone of the first is normal, and the 
affine cone of the second is not.

In Exercise II.5.14, \cite{_Hartshorne:AG_},
Hartshorne gives a criterion for projective normality:

\hfill

\proposition\label{proj_normality_Proposition_}
A projective variety $X\subset \C P^n$ is projectively normal if and
only if $X$ is normal, and the restriction map 
\[ H^0(\C P^n, \calo(i)) \arrow H^0\left(X, \calo(i)\restrict X\right)\]
is surjective for all $i \geq 0$. 
\endproof 

\hfill

\remark
This proposition can be applied to algebraic cones,
because the quotient singularities are normal 
(\ref{_quotient_normal_Claim_}),
and the orbifolds have only quotient singularities.

\hfill

\example
\label{_non_unique_non_normal_Example_}
Let $X\subset \C P^n$ be a projective manifold,
which is not contained in a projective subspace
$\C P^k \subset \C P^n$ of smaller dimension.
Consider a linear projection $p:\; \C P^n \arrow \C P^{n-1}$
centered in a point $z\notin X$. Then the restriction $p\restrict X$ is holomorphic.
Assume that $p:\; X \arrow \C P^{n-1}$ is also injective
(this can be always achieved if $2\dim X < n-1$ for 
an appropriately general choice of $z$; indeed, $p$ is not injective
if and only if the center $z$ does not belong to the secant 
variety\footnote{A secant variety of $X\subset \C P^n$ is the Zariski closure of the union
 of all
lines $\C P^n \subset \C P^n$ intersecting $X$ in at least two points.} of $X$,
which has dimension $\leq 2\dim X +1$). Consider the natural map of affine cones
$u:\; C(X) \arrow C(p(X))$. Outside of the origin $c$, this map is biholomorphic,
hence it is birational and finite. By definition, 
a variety $Z$ is normal
if any bimeromorphic finite map $Z_1 \arrow Z$ is an isomorphism.
Were $C(p(X))$ normal, this would imply
that $u$ is an isomorphism. However, the Zariski tangent space
$T_c C(X)$ is $n+1$-dimensional, because $C(X)\subset C(\C P^n)$ generates 
the vector space $C(\C P^n)= \C^{n+1}$ (otherwise $X$ would have been contained
in a smaller dimensional projective subspace). On the other hand,
$\dim T_c C(p(X)) \leq n$, hence $u:\; C(X) \arrow C(p(X))$
is not an isomorphism.\footnote{We are grateful to Yu. Prokhorov for
  this example.}

\hfill

This example makes sense, if one considers the following statement.

\hfill

\claim\label{_closed_cone_normal_from_open_Claim_}
Let $\tilde M$ be an open algebraic cone, and
${\goth S}$ the set of all closed algebraic cones
$\tilde M_c$ obtained by adding the origin to $\tilde M$.
Then there exists only one closed algebraic cone $Z\in {\goth S}$ which
is normal, and for any other $Z' \in {\goth S}$, 
its normalization is $Z$. Moreover, the normalization
map $Z \arrow Z'$ is a homeomorphism. 

\hfill

\proof
Let $Z$ be the Stein completion of $\tilde M$;
by definition, it is normal. Since all holomorphic functions
on $\tilde M$ can be extended to $Z$, this variety is equipped with
a finite, bijective, bimeromorphic map to any 
closed algbraic cone $Z'$ associated with $\tilde M$.
Indeed, the holomorphic functions on $Z'$ are
holomorphic on $\tilde M$, hence they can be
extended to its Stein completion.
This implies that $Z$ is the normalization of all other 
$Z' \in {\goth S}$.
\endproof

\hfill

{\bf Acknowledgements:}
We are grateful to Yuri Prokhorov for 
\ref{_non_unique_non_normal_Example_} and
an important
discussion about the projective normality,
and to Andrey Soldatenkov for his insightful comments
about the algebraic structures.

\hfill

{\scriptsize

\noindent {\sc Liviu Ornea\\
University of Bucharest, Faculty of Mathematics and Informatics, \\14
Academiei str., 70109 Bucharest, Romania}, and:\\
{\sc Institute of Mathematics ``Simion Stoilow" of the Romanian
Academy,\\
21, Calea Grivitei Str.
010702-Bucharest, Romania\\
{\tt lornea@fmi.unibuc.ro,   liviu.ornea@imar.ro}}

\hfill

\noindent {\sc Misha Verbitsky\\
Instituto Nacional de Matem\'atica Pura e
Aplicada (IMPA) \\ Estrada Dona Castorina, 110\\
Jardim Bot\^anico, CEP 22460-320\\
Rio de Janeiro, RJ - Brasil,} also:\\
{\sc Laboratory of Algebraic Geometry, \\
Faculty of Mathematics, National Research University 
Higher School of Economics,
6 Usacheva Str. Moscow, Russia\\
{\tt verbit@verbit.ru, verbit@impa.br}}

\hfill

}

\end{document}